\documentclass[11.5pt]{amsart}
\usepackage{graphicx}
\usepackage{amssymb}
\usepackage{amsmath}
\usepackage{amsthm}
\usepackage{amscd}

\newtheorem{thm}{Theorem}[section]
\newtheorem{sub}[thm]{}
\newtheorem{cor}[thm]{Corollary}
\newtheorem{lem}[thm]{Lemma}
\newtheorem{exm}[thm]{Example}

\newtheorem{prop}[thm]{Proposition}
\theoremstyle{definition}

\theoremstyle{remark}
\newtheorem{rem}[thm]{Remark}
\numberwithin{equation}{section}

\begin{document}
\title[Representations and Cocycle Twists  of Color Lie Algebras ]
{Representations and Cocycle Twists  of  Color Lie Algebras}
\author[Xiao-Wu Chen, Sergei D. Silvestrov and Fred van Oystaeyen
] {X.-W. Chen$^{a, c}$, S. D. Silvestrov $^b$ and F. van Oystaeyen$^{c,*}$ }

\thanks{$^*$ The corresponding author}
\thanks{Supported by
AsiaLink project  ``Algebras and Representations in China and
Europe$"$ ASI/B7-301/98/679-11, and the ESF Scientific Programme
``NOG'', and National  Natural Science Foundation of China No.
10271113.}

\maketitle

\begin{center}
$^A$Department of Mathematics \\University of Science and
Technology of China \\Hefei 230026, Anhui, P. R. China
\end{center}

\vskip5pt
\begin{center}
$^B$Centre for Mathematical Sciences, Department of Mathematics\\
Lund Institute of Technology,  Lund University\\
Box 118, SE-221 00 Lund, Sweden
\end{center}
\vskip 5pt
\begin{center}
 $^C$Department of Mathematics and Computer Science\\
University of Antwerp, Middelheimlaan 1\\
 B-2020 Antwerp, Belgium
\end{center}

\vskip 5pt
\begin{center}
xwchen$\symbol{64}$mail.ustc.edu.cn\\
sergei.silvestrov$\symbol{64}$math.lth.se\\
fred.vanoystaeyen$\symbol{64}$ua.ac.be
\end{center}
\begin{abstract}
We study relations between finite-dimensional representations of color Lie
algebras and their cocycle twists. Main tools are the universal enveloping
algebras and their FCR-properties (finite-dimensional representations are
completely reducible.) Cocycle twist preserves the FCR-property. As an
application, we compute all finite dimensional representations (up to
isomorphism) of the color Lie algebra $sl_2^c$.
\end{abstract}

\vskip 20pt

\section{{\bf Introduction}}
The generalizations of Lie algebras and Lie superalgebras introduced in [RW]
and systematically studied in  [Sc1], crossed over from physics to abstract
algebra. Nowadays they are known as color Lie algebras. \par \vskip 10pt

Throughout $K$ is  a field of characteristic zero and $L$ is a color Lie
algebra over $K$.  Two associative $K$-algebras may be associated to $L$, the
universal enveloping algebra $\mathcal{U}(L)$ and the augmented  enveloping
algebra $\widetilde{\mathcal{U}}(L)$ (see  2.7).\par
The category of representations of $L$ is equivalent to the category of
$\mathcal{U}(L)$-modules. Together with $\mathcal{U}(L)$ the augmented
enveloping algebra $\widetilde{\mathcal{U}}(L)$ appears naturally, the latter
is a Hopf algebra containing $\mathcal{U}(L)$ as a subalgebra.
\par\par \vskip 10pt
An associative $K$-algebra $A$ is an FCR-algebra, following [KS] and  [KSW], if
every finite-dimensional representation is completely reducible and the
intersection of its annihilators of all the finite-dimensional representations
is zero. Classical examples of FCR-algebras are finite-dimensional semisimple
algebras, the universal enveloping algebra $\mathcal{U}({\bf g})$ of a
finite-dimensional semisimple  Lie algebra ${\bf g}$, the quantum enveloping
algebras $\mathcal{U}_q({\bf g})$ with  $q$  not a root of unity, and the
universal enveloping algebra $\mathcal{U}(osp(1, 2r))$ of the orthosymplectic
Lie superalgebra $osp(1,2r)$ (references for these facts are in the
Introduction and Theorem 2.13 of [KiS]).
\par \vskip 10pt
It is then natural to ask which  color Lie algebras $L$ are such that
$\mathcal{U}(L)$ and $\widetilde{\mathcal{U}}(L)$ are FCR-algebras. We provide
a partial answer to  this question in this note. We restrict attention to
finite-dimensional color Lie algebra $L$ graded by a finite abelian group $G$.
The main observations are
\par \vskip 10pt
\begin{thm} Let $L$ be a finite-dimensional $G$-graded $\varepsilon$-Lie
algebra for a finite abelian group $G$. Then we have  \vskip 5pt
\noindent {\rm (1)}. If $L^c$ is a cocycle twist of $L$, then $\mathcal{U}(L^c)$
is an FCR-algebra if and only if $\mathcal{U}(L)$ is an FCR-algebra. \\
{\rm (2)}. $\mathcal{U}(L)$ is an FCR-algebra if and only if
$\widetilde{\mathcal{U}}(L)$ is an FCR-algebra.
\end{thm}
\par \vskip 10pt
It follows immediately from Theorem 1.1 (1)  that cocycle twists of a
finite-dimensional semisimple Lie algebra have FCR-property for their universal
enveloping algebras and augmented enveloping algebras. \par \vskip 10pt

The paper is arranged as follows: Section 2 contains some preliminaries on color
Lie algebras; Section 3 is devoted to the notion of cocycle twists of color
Lie algebras, which turns out to coincide with the cocycle twists of graded
associative algebras (Theorem 3.5). In Section 4  we prove some results on
FCR-algebras (e.g. Theorem 4.2), and apply these to show Theorem 1.1. Section 5
contains the computation of all finite-dimensional representations  of a
specific color Lie algebra $sl_2^c$, which is a cocycle twist of  the simple
Lie algebra $sl(2, K)$ for $K$ an algebraically closed field(see Example 3.3
and Theorem 5.2).

\vskip 20pt

\section{{\bf Preliminaries}}

Throughout we will work over a field $K$ of characteristic zero. All unadorned
tensor products are over $K$. By a  noetherian algebra we mean a two-sided
noetherian $K$-algebra.
\par \vskip 5pt

\sub{} \rm Let us recall some basic facts on color Lie algebras ([Sc1]).
Suppose that  $G$ is  an abelian group with the identity element $e$, and
$\varepsilon: G \times G \longrightarrow K^\times$ is an antisymmetric
bi-character, i.e.,
\begin{align*}
&\varepsilon(g, h) \varepsilon(h, g)=1,\\
&\varepsilon(g, hk)=\varepsilon(g, h) \varepsilon(g, k),\\
&\varepsilon(gh, k)=\varepsilon(g, k) \varepsilon(h, k),
\end{align*}
for all $g, h, k \in G$, where  $K^\times=K\setminus \{0\}$ is the
multiplicative group of units in $K$. In particular, we note that
$\varepsilon(g, g)=1$ or $-1$, $g \in G$.\par \vskip 5pt

By a $G$-graded $\varepsilon$-Lie algebra $L$ (or a color Lie algebra), we mean
that $L=\oplus_{g \in G}L_g$ is a graded space over $K$, equipped with a
bilinear multiplication
\begin{align*}
\langle -,- \rangle: L \times L \longrightarrow L
\end{align*}
such that
\begin{align*}
&\langle L_g, L_h \rangle  \subseteq L_{gh}\quad
(\mbox{gradation condition}); \\
&\langle x, y \rangle =-\varepsilon(g, h)\langle y, x\rangle \quad
(\mbox{color symmetry}); \\
&\varepsilon(k, g)\langle x,\langle  y , z \rangle \rangle +\varepsilon(g, h)
\langle y, \langle  z, x\rangle \rangle + \varepsilon(h, k)\langle z,
\langle  x , y\rangle \rangle =0 \quad  (\mbox{color Jacobi identity}),
\end{align*}
for all $x \in L_g$, $y \in L_h$ and $z \in L_k$, where  $g, h, k \in G$.
\par \vskip 10pt
For example, if $G=\mathbb{Z}_2$ and $\varepsilon(\alpha, \beta)=(-1)^{\alpha
\beta}$ for all $\alpha, \beta \in G$, then $G$-graded $\varepsilon$-Lie
algebras are exactly Lie superalgebras ([Sc2]); if $\varepsilon$ is trivial,
i.e., $\varepsilon(g, h)=1$ for all $g, h \in G$, then $G$-graded
$\varepsilon$-Lie algebras are known as $G$-graded Lie algebras. \par \vskip 5pt
An element $x \in L$ is said to be homogeneous of degree $g$ if $x \in L_g$,
and  in this case,  we write $|x|=g$. So the gradation condition implies
that $|\langle x, y\rangle|=|x|\cdot|y|$ for homogeneous elements $x$ and $y$
in $L$ (here, the dot denotes the multiplication in $G$). \par \vskip 10pt

Recall that a representation $(V, \rho)$ of a $G$-graded $\varepsilon$-Lie
algebra $L$ is just a linear map $\rho: L \longrightarrow {\rm End}_K(V)$ such
that
\begin{align*}
\rho(\langle x, y \rangle )=\rho(x) \rho(y) -\varepsilon(g, h) \rho(y)\rho(x),
\end{align*}
for all $x \in L_g$, $y \in L_h$ and $g, h \in G$. \par
A representation $(V, \rho)$ is said to be graded, if $V=\oplus_{g \in G} V_g$
is a $G$-graded space such that
\begin{align*}
\rho(x) (V_g) \subseteq V_{|x|g}
\end{align*}
for all the homogeneous elements $x \in L$ and $g \in G$. In the present paper,
we will not restrict ourselves to graded representations.
\par \vskip 10pt
\sub{$\varepsilon$-symmetric algebras $S_{\varepsilon} (X)$}\vskip 10pt
\rm  Let $X$ be a finite set. The set $X$ is said to be labelled by an abelian
group $G$, if there is  a map $| \cdot |: X \longrightarrow G$.  Let
$\varepsilon : G \times G \longrightarrow K^\times$ be any bilinear form.
Then the $\varepsilon$-symmetric algebra $S_\varepsilon (X)$ on $X$
([Section 4.B, Sc1]) is an associative $K$-algebra  generated by elements in
$X$ subject to relations
\begin{align*}
xy= \varepsilon(|x|, |y|) yx, \quad \forall x, y \in X.
\end{align*}
\vskip 5pt
A well-known fact is
\begin{lem}
The $\varepsilon$-symmetric  algebra $S_\varepsilon (X)$ is noetherian.
\end{lem}
\noindent{\bf Proof} \quad Use induction on the order $|X|$ of the set $X$.
If $|X|=1$, say $X=\{x\}$, then  $S_\varepsilon(X)=k[x]$ if $\varepsilon(|x|,
|x|)=1$;  $S_{\varepsilon}(X)=k[x]/{(x^2)}$ if $\varepsilon(|x|, |x|)\neq 1$,
both of which  are  clearly noetherian. \par
Assume that $x \in X$ and   $X'=X \setminus \{x\}$. By induction
$A_\varepsilon(X')$ is noetherian. Let $\sigma_x $ be the automorphism of
$A_\varepsilon(X')$ given by
\begin{align*}
\sigma_x (y)= \varepsilon(|x|, |y|) y, \quad \forall y \in X'.
\end{align*}
Denote by  $A_\varepsilon(X')[t; \sigma_x]$ the Ore extension of
$A_\varepsilon(X')$ ([MB], p.15) with respect to the automorphism $\sigma_x$,
where $t$ is a variable.  By [Theorem 1.2.9, MB],
$A_\varepsilon(X')[t; \sigma_x]$ is noetherian. \par

Note that  there is an epimorphism of algebras
\begin{align*}
\pi: A_\varepsilon(X')[t; \sigma_x] \longrightarrow S_\varepsilon(X)
\end{align*}
such that $\pi(t)=x$ and $\pi(y)=y$,  $y \in X'$. It follows that
$S_\varepsilon (X)$ is noetherian. \hfill $\blacksquare$
\vskip 10pt
\sub{Universal enveloping algebras $\mathcal{U}(L)$}\vskip 10pt
\rm
Let $L$ be a $G$-graded $\varepsilon$-Lie algebra as in 2.1. Recall that the
universal enveloping algebra $\mathcal{U}(L)$ of $L$ is defined to be
\begin{align*}
\mathcal{U}(L)=\mathcal{T}(L)/{\mathcal{J}(L)}
\end{align*}
where $\mathcal{T}(L)$ is the tensor algebra of $L$ and $\mathcal{J}(L)$ is the
two-sided ideal generated by $x \otimes y -\varepsilon(g,h) y\otimes x - \langle
x, y  \rangle$ for all $x \in L_g$, $y \in L_h$ and $g, h \in G$.\par \vskip 5pt
Note that if $L$ is trivial, i.e., $\langle L, L \rangle=0$, then
$\mathcal{U}(L) \simeq S_{\varepsilon}(X)$, where the set  $X$ is a homogeneous
basis of $L$  labelled by the degrees of these elements in $L$. \par \vskip 5pt

Write the canonical map $i_L: L \longrightarrow \mathcal{U}(L)$.  Note that
$\mathcal{U}(L)$ is a $G$-graded associative algebra such that $i_L$ is a
graded map of degree $e$. For theory on $G$-graded rings, we refer to the book
[NV2].\par \vskip 5pt

Obviously here is an equivelence between the category of representations of the
color  Lie algebra $L$ and the category of $\mathcal{U}(L)$-modules. Also
graded representations of $L$ correspond exactly to graded
$\mathcal{U}(L)$-modules (note again that $\mathcal{U}(L)$ is a $G$-graded
algebra).  So the question we consider  in the introduction can be re-stated
as: which color Lie algebras have enough finite-dimensional representations and
all of these being completely reducible. \par \vskip 10pt

A remarkable fact is  that the  Poincar\'e-Birkhoff-Witt theorem holds for the
algebra $\mathcal{U}(L)$ of any color Lie algebra $L$, see [Theorem 1, Sc1] or
[Section 3, Hu]. For later use, we quote it as
\par \vskip 10pt
\begin{prop}
Let $L$ be a $G$-graded $\varepsilon$-Lie algebra with the universal enveloping
algebra $\mathcal{U}(L)$, then the canonical map $i_L$ is injective. Moreover,
if $\{x_i\}_{i \in \Lambda}$ is a homogeneous basis of $L$, where $\Lambda$ is
a well-ordered set, then the set of  ordered  monomials
$x_{i_1}x_{i_2}\cdots x_{i_n} $ is a basis of $\mathcal{U}(L)$, where
$i_j \leq i_{j+1}$ and $i_j < i_{j+1}$ if $\varepsilon(|x_{i_j}|, |x_{i_j}|)
=- 1$,  $1 \leq j \leq n$ and  $n = 0, 1, \cdots$.
\end{prop}
\par \vskip 10pt
\noindent A direct consequence is
\begin{cor} If $L$ is a finite-dimensional $G$-graded $\varepsilon$-Lie algebra,
then $\mathcal{U}(L)$ is a noetherian algebra.
\end{cor}
\noindent{\bf Proof} \quad  By [Section 4.C, Sc1], we see that $\mathcal{U}(L)$
is a positively filtered algebra with its  associated graded algebra $gr
\mathcal{U}(L) \simeq S_\varepsilon (X) $, where $X$ is a set of homogeneous
basis of $L$ and is labelled by the degrees of the basis elements in $L$. It
follows from Lemma 2.3 and [Theorem 1.6.9, MB] that $\mathcal{U}(L)$ is
noetherian. \hfill $\blacksquare$
\newpage
\sub{Augmented enveloping algebras $\widetilde{\mathcal{U}}(L)$} \vskip 5pt
\rm
Let us recall the notion of an augmented  enveloping algebra
$\widetilde{\mathcal{U}}(L)$ of the $G$-graded $\varepsilon$-Lie algebra $L$:
as a vector space $\widetilde{\mathcal{U}}(L):=KG \otimes \mathcal{U}(L)$,
where $KG$ is the group algebra of $G$. The multiplication is given as
\begin{align*}
(g \otimes x) \cdot (h \otimes y):=\varepsilon(|x|, h) gh \otimes xy
\end{align*}
where $g, h\in G$ and  $x, y \in \mathcal{U}(L)$ are homogeneous. (Note that
$\mathcal{U}(L)=\oplus_{g \in G} \mathcal{U}(L)_g$ is a $G$-graded algebra).
One sees that $\widetilde{\mathcal{U}}(L)$ contains $\mathcal{U}(L)$ and $KG$
as subalgebras.\par \vskip 5pt
The algebra $\widetilde{\mathcal{U}}(L)$ is a Hopf algebra with coalgebra
structure maps $\Delta$, $\epsilon$ and antipode $S$ given by
\begin{align*}
\Delta(g)=g \otimes g,& \quad \Delta(x)=1 \otimes x + x \otimes h\\
\epsilon(g)=1, &\quad \epsilon(x)=0\\
S(g)=g^{-1}, & \quad S(x)=-xh^{-1},
\end{align*}
where $g \in G$ and $x \in L_h$. (See [Hu] or [Kha].)\par\vskip 5pt
We have  \par
\begin{prop}
Let $L$ be a $G$-graded $\varepsilon$-Lie algebra. If $L$ is finite-dimensional
and $G$ is a finite group, then $\widetilde{\mathcal{U}}(L)$ is a  noetherian
Hopf algebra.
\end{prop}
\noindent{\bf Proof} \quad Since $\widetilde{\mathcal{U}}(L)$ is a free left
and right $\mathcal{U}(L)$-module of finite rank, and $\mathcal{U}(L)$ is
noetherian by Corollary 2.6, it follows immediately that
$\widetilde{\mathcal{U}}(L)$ is a  noetherian left and right
$\mathcal{U}(L)$-module, so it is a noetherian algebra. \hfill $\blacksquare$
\par
\vskip 20pt
\section{{\bf Cocycle Twist}}
\vskip 5pt
In this section, we study cocycle twists of color Lie algebras, which turn out
to correspond to  the cocycle twists of $G$-graded  associative algebras, see
Theorem 3.5. We also offer an explicit example of a color Lie algebra $sl_2^c$.
\par \vskip 5pt
\sub{} \rm  Let us recall the notion of a cocycle twist of a color Lie algebra
as  given in  [Sc1] (also see [RW]).
Let  $c: G \times G \longrightarrow K^\times$ be a cocycle on the abelian group
$G$ with trivial action of $G$ on $K$, i.e.,
\begin{align}
c(g, h) c(gh, k)= c(h, k) c(g, hk),
\end{align}
for all $g, h, k \in G$. Define
\begin{align*}
B_c(g, h)=\frac{c(g, h)}{c(h, g)}, \quad \forall g, h \in G
\end{align*}
which is a bilinear form of $G$, see [Kar]. \par \vskip 10pt

Let $L=\oplus_{g\in G} L_g$ be a $G$-graded $\varepsilon$-Lie algebra. Fix a
triple $(c, \phi, \varepsilon')$,  where $c$ is a cocycle of $G$,  $\phi: G
\longrightarrow G'$ is a morphism of abelian groups, and $\varepsilon'$ is
bi-character of $G'$ such that
\begin{equation}
\varepsilon(g,h)c(g, h)=c(h,g) \varepsilon'(\phi(g), \phi(h)).
\end{equation}
\vskip 10 pt

The cocycle twist ([Section 6, Sc1]) of $L$ with respect to the triple $(c,\phi,
\varepsilon')$, denoted by $L^c$,  is defined as follows:  \\
(1).\quad $L^c=L$ as $K$-spaces;  $L^c$ is $G'$-graded by setting
$L^c_{g'}=\oplus_{\phi(g)=g'}L_g$ for all $g' \in G'$. \\
(2). \quad Define the bracket $\langle-,- \rangle^c$ on $L^c$ as
\begin{align*}
\langle x, y \rangle^c := c(g, h) \langle x, y \rangle,
\end{align*}
for all  $x \in L_g \subseteq L^c_{\phi(g)}$,  $ y \in L_h \subseteq
L^c_{\phi(h)}$ and  $g, h \in G$. \par \vskip 10pt
\begin{lem}
With notation as above,  $L^c$ is a $G'$-graded $\varepsilon'$-Lie algebra.
\end{lem}
\noindent {\bf Proof} \quad We have to check the three axioms in the definition
of  color Lie algebras. The gradation condition of $L^c$ holds (since $\phi$ is
a morphism);  and the color symmetry follows  from the equation (3.2).\par
What is left is the  color Jacobi identity, i.e.,
\begin{align*}
\varepsilon'(\phi(k), \phi(g) ) \langle x, \langle y, z\rangle^c \rangle^c +
\varepsilon'(\phi(g), \phi(h))\langle y, \langle z, x\rangle^c \rangle^c+
\varepsilon'(\phi(h), \phi(k))\langle z, \langle  x , y\rangle^c \rangle^c =0
\end{align*}
for all $x \in L_g \subseteq L^c_{\phi(g)}$, $y\in L_h \subseteq L^c_{\phi(h)}$,
$z \in L_k \subseteq L^c_{\phi(k)}$ and  $g, h, k \in G$.\par \vskip 5pt
By the definition of $\langle -, - \rangle^c$ and the color Jacobi identity of
$L$ itself, it suffices to show
\begin{align*}
\frac{\varepsilon'(\phi(k), \phi(g))}{\varepsilon(k, g)} c(h, k) c(g, hk)&=
\frac{\varepsilon'(\phi(g), \phi(h))}{\varepsilon(g, h)} c(k, g) c(h, kg)\\
&=\frac{\varepsilon'(\phi(h), \phi(k))}{\varepsilon(h, k)} c(g, h) c(k, gh),
\end{align*}
for all $g, h , k \in G$. \par \vskip 5pt
\noindent By (3.2), we have
\begin{align*}
\frac{\varepsilon'(\phi(k),\phi(g))}{\varepsilon(k,g)}=\frac{c(k, g)}{c(g, k)}
\quad \mbox{and} \quad  \frac{\varepsilon'(\phi(g), \phi(h))}{\varepsilon(g, h)}
=\frac{c(g, h)}{c(h,g)}.
\end{align*}
So the first equality follows from
\begin{align*}
\frac{c(k, g)}{c(g,k)} c(h,k)c(g, hk)= \frac{c(g, h)}{c(h, g)} c(k, g) c(h, kg).
\end{align*}
Now apply $c(k, g) c(h, kg)=c(h, k) c(hk, g)$ at the right side of the above
equation, we obtain that the first equality follows exactly from
\begin{align*}
\frac{c(h, g)}{c(g, h)} \cdot \frac{c(k, g)}{c(g , k)} = \frac{c(hk, g)}
{c(g, hk)},
\end{align*}
which is just the bilinearity of $B_c$ (see 3.1). Similarly, one can show the
second equality. This completes the proof. \hfill $\blacksquare$
\vskip 10pt
Let us consider an example of cocycle twist, which will be studied further in
Section 5.  For more examples of cocycle twists, see [RW] and [Sc1].
\par \vskip 5pt
\begin{exm}
\rm Assume that $K$ is algebraically closed. $sl_2=sl(2, K)$ is the
three-dimensional simple Lie algebra over $K$, with stardard basis $\{e,h,f\}$
such that
\begin{align*}
[h, e]=2e, \quad [e, f]=h, \quad [h, f]=-2f,
\end{align*}
where $[-, -]$ denotes the Lie bracket of $sl_2$. \par
Put $a_1=\frac{i}{2} (e-f)$, $a_2=-\frac{1}{2}(e+f)$ and $a_3=\frac{i}{2} h$
(where $i^2=-1$) in $sl_2$. So
\begin{align*}
[a_1, a_2]=-a_3, \quad [a_2, a_3]=-a_1, \quad [a_3,a_1]=a_2.
\end{align*}
Assume that $G=\mathbb{Z}_2 \times \mathbb{Z}_2$.
So $sl_2=\oplus_{g \in G}X_g$ is a $G$-graded Lie algbera, with
\begin{align*}
X_{(0, 0)}=0,\quad X_{(1, 0)}=Ka_1,\quad X_{(0, 1)}=Ka_2,\quad X_{(1, 1)}=Ka_3.
\end{align*}
Let $c$ be a cocycle of $G$ given by
\begin{align*}
c((\alpha_1, \alpha_2), (\beta_1, \beta_2))=(-1)^{\alpha_1 \beta_2}
\end{align*}
for all $\alpha_1, \alpha_2, \beta_1, \beta_2 \in \mathbb{Z}_2$.  \par
Take $G'=G$, $\phi=Id$ and the bi-character $\varepsilon'$ of $G'=G$  to be
\begin{align*}
\varepsilon'((\alpha_1, \alpha_2), (\beta_1, \beta_2)):= (-1)^{\alpha_1
\beta_2-\alpha_2 \beta_1}
\end{align*}
for all $\alpha_1, \alpha_2, \beta_1, \beta_2 \in \mathbb{Z}_2$. Clearly the
equation (3.2) holds (here the $\varepsilon$ is trivial). \par \vskip 5pt
Denote by $sl_2^c$  the cocycle twist of $sl_2$ with respect to the triple
$(c,Id, \varepsilon')$. By Lemma 3.2, $sl_2^c$ is a $G$-graded
$\varepsilon'$-Lie algebra. The algebra $sl_2^c$ has a homogeneous basis
$\{a_1, a_2, a_3\}$, with degrees given by
\begin{align*}
|a_1|=(1, 0), \quad |a_2|=(0,1), \quad|a_3|=(1,1).
\end{align*}
By definition of the cocycle twist, we obtain that the bracket
$\langle-, - \rangle$ in $sl_2^c$ is given by:
\begin{align*}
\langle a_1, a_2 \rangle =a_3, \quad \langle a_2, a_3 \rangle =a_1, \quad
\langle a_3, a_1 \rangle =a_2.
\end{align*}
In other words, the universal enveloping algebra $\mathcal{U}(sl_2^c)$ is
generated by $a_1, a_2, a_3$ with relations
\begin{align*}
a_1 a_2 + a_2 a_1 &= a_3,  \\
a_2 a_3 + a_3 a_2 &= a_1, \\
a_3 a_1 + a_1 a_3 &= a_2,
\end{align*}
which is a $G$-graded associative algebra,  with $|a_1|=(1, 0)$, $|a_2|=(0,1)$
and $|a_3|=(1,1)$. Note that this algebra and its representations appear in
[Appendix, RW], [OS] and [GP].
\end{exm}
\vskip 10pt
\sub{} \rm  We will show that cocycle twists of color Lie algebras in 3.1
coincide with the well-known notion of cocycle twists (e.g., see [AST]) of
associative algebras, see Theorem 3.5. \par \vskip 10pt

Let $G$ be an abelian group and  $A=\oplus_{g \in G} A_g$ be a $G$-graded
associative $K$-algebra. For any cocycle $c$ of $G$, the cocycle twist of $A$
with respect to $c$, denoted by $A^c$, is defined as: $A^c:=A$ as $G$-graded
spaces; the multiplication $``\circ"$ on $A^c$ is given by
\begin{align*}
x^c \circ y^c :=c(g, h) (xy)^c, \quad x^c \in A^c_g, y \in A^c_h,
\end{align*}
where $a^c$ denotes the element in $A^c$ which corresponds to $a \in A$.
\par \vskip 5pt

For example, if $A=KG$ is the group algebra with the natural grading by $G$,
then $KG^c$ is exactly the twisted group ring with respect to the cocycle $c$
([NV2], p.12).

\par \vskip 10pt

The main result of this section is
\begin{thm}
Using notation as in {\rm 3.1} and the above, we have an isomorphism
$$\mathcal{U}(L)^c \simeq \mathcal{U}(L^c),$$
where $\mathcal{U}(L)$ is considered as a $G$-graded algebra and
$\mathcal{U}(L)^c$ its cocycle twist with respect to  $c$.
\end{thm}

\noindent{\bf Proof}\quad Identify $L$ as a subspace of $\mathcal{U}(L)$,
hence a subspace of $\mathcal{U}(L)^c$. We claim that there is an algebra
morphism
\begin{align*}
\Theta: \mathcal{U}(L^c) \longrightarrow \mathcal{U}(L)^c
\end{align*}
such that $\Theta(x)=x^c$ for all $x \in L^c$. (Note that $L^c=L$ as vector
spaces.)\par

By the definition of $\mathcal{U}(L^c)$, it suffices to check that
\begin{align*}
\Theta(x)\circ\Theta(y)-\varepsilon'(\phi(g), \phi(h))\Theta(y)\circ
\Theta(x)=\Theta(\langle x, y \rangle^c ),
\end{align*}
for all $x \in L_g \subseteq L^c_{\phi(g)}$ and $y \in L_h \subseteq
L^c_{\phi(h)}$, where $``\circ"$ denotes the multiplication of
$\mathcal{U}(L)^c$ and $\langle-, -\rangle^c$ the bracket of $L^c$.\par
In fact,
\begin{align*}
&\Theta(x)\circ\Theta(y)-\varepsilon'(\phi(g),\phi(h))\Theta(y)\circ \Theta(x)\\
&= x^c \circ y^c -\varepsilon'(\phi(g), \phi(h)) y^c \circ x^c\\
&=c(g, h)(xy)^c - \varepsilon'(\phi(g), \phi(h)) c(h, g)(yx)^c\\
&=c(g,h) (xy)^c- c(g, h) \varepsilon(g, h) (yx)^c\\
&=c(g,h) (\langle x, y \rangle )^c\\
&=\Theta(\langle x, y \rangle ^c).
\end{align*}
(The fourth equality uses the fact that $xy-\varepsilon(g, h) yx=\langle x, y
\rangle$ in $\mathcal{U}(L)$; the last one follows from the definition
$\langle x, y \rangle^c=c(g, h) \langle x, y \rangle$.) So we have shown that
the algebra morphism $\Theta$ is well-defined.\par \vskip5pt

Clearly $\Theta$ is surjective, since $\mathcal{U}(L)^c$ is generated by the
image $\Theta(L^c)$ as an associative algebra.\par
To show that  $\Theta$ is injective, let $\{x_i\}_{i \in \Lambda}$ be a set of
homogeneous basis of $L$ (hence of $L^c$) with $\Lambda$ well-ordered, then by
Proposition 2.5,  the  monomials $x_{i_1}x_{i_2} \cdots x_{i_n}$,
$i_j \leq i_{j+1}$ and $i_j < i_{j+1}$ whenever $\varepsilon'(\phi(g_j),
\phi(g_j))=- 1$, where $x_j \in L_{g_j} \subseteq L^c_{\phi(g_j)}$ ($1 \leq j
\leq n$), forms a basis of $\mathcal{U}(L^c)$. \par
Note that $\varepsilon(g, g)=\varepsilon'(\phi(g), \phi(g))$ by equation (3.2).
In particular, $\varepsilon(g, g)=1$ if and only if $\varepsilon'(\phi(g),
\phi(g))=1$. So again by Proposition 2.5,  the set of the images of these
monomials, i.e.,  $\{\Theta(x_{i_1}x_{i_2} \cdots x_{i_n})\}$,  is linearly
independent in $\mathcal{U}(L)^c$, it follows that $\Theta$ is injective. This
completes the proof. \hfill $\blacksquare$

 \vskip 20pt

\section{{\bf FCR-algebras and Proof of Theorem 1.1}}
\vskip 10pt
In this section, we will prove  several  results on FCR-algebras, from which Theorem 1.1 follows.\par \vskip 5pt

   \sub{} \rm  Let    $A_e$ be a $K$-algebra and  $G$ a group  acting on $A_e$ as algebra automorphisms, and let $c: G \times G \longrightarrow K^\times$ be a cocycle. Recall that the  crossed product ([NV2], p.11)  $A_e *_c G$ is just a free $A_e$-module with basis $\bar{G}$, where $\overline{G}=\{\overline{g}| g \in G\}$ is a copy of $G$ as a set, and its multiplication is given by
   \begin{center}
   $(a * \overline{g}) \cdot (b * \overline{h}):= c(g, h)a (g.b)* \overline{gh}$,
   \end{center}
   for all  $g, h \in G$ and $a, b \in A_e$. (We denote by ``.''  the $G$-action on $A_e$.)\par\vskip 5pt

   Recall  a notation (e.g., see [NV1], p77): let $A$ be an algebra and $B\subseteq A$ a subalgebra, and  let $M$ be a left $B$-module. Then $A \otimes_B M$ becomes a left $A$-module defined by
    \begin{align*}
    a. a'\otimes_B m:= aa' \otimes_B m
    \end{align*}
    for all $a, a' \in A$ and $m \in M$, where $a \otimes_B m$ is viewed as an element of $A\otimes_B M$. Moreover, if $M$ is a finite-dimensional $B$-module and  $A$ is a finitely-generated right $B$-module, then $A\otimes_B M$ is also  finite-dimensional.\par \vskip 5pt

    We need  the following result.\par

   \begin{thm}
   Let $G$ be a finite abelian group  and $A_e$ a noetherian $K$-algebra. Suppose that  $A=A_e *_c G$ is a crossed product as above, then $A$ is an FCR-algebra if and only if $A_e$ is an FCR-algebra.
   \end{thm}

   \noindent{\bf Proof}\quad We may assume that $c(g, e)=c(e, g)=1$ for all $g \in G$. (Otherwise, replace $c$ by $\frac{c}{c(e, e)}$, note that for a cocycle $c$,  $c(g, e)=c(e, h)$ for all $g, h \in G$, see the Appendix.) Identify $A_e$ as a subalgebra of $A=A_e *_c G$ by sending $a$ to $a * \overline{e}$ for all $a \in A_e$, and view $\overline{G}$ as a subset of $A$ by identifying $\overline{g}$ with $1 * \overline{g}$ for all $g \in G$.  \par \vskip 5pt

   For the ``if" part, assume that $A_e$ is an FCR-algebra. Let $W$ be a finite-dimensional $A$-module with a submodule $V$. So $W$ is also an $A_e$-module with the $A_e$-submodule $V$. Since $A_e$ is an FCR-algebra, in particular, every finite-dimensional $A$-module is completely reducible. So $V$ is a direct summand of $W$ as $A_e$-modules. In other words, there exists an $A_e$-module morphism
   \begin{align*}
   p_e: W \longrightarrow V
   \end{align*}
   such that $p_e|_V=Id_V$. \par
   Define $p:W \longrightarrow V$ as
   \begin{align*}
   p(w):=\frac{1}{|G|} \sum_{g \in G} \frac{1}{c(g, g^{-1})} \overline{g}p_e(\overline{g^{-1}}w), \quad \forall w \in W.
   \end{align*}
   We claim that $p$ is an $A$-module morphism. \par \vskip 5pt

   Note that $A$ is generated by $A_e$ and the set $\overline{G}$ as an associative algebra. So it suffices to show $p(aw)=ap(w)$ and $p(\overline{h}w)=\overline{h} p(w)$, for all $a \in A_e$, $h \in G$ and $w \in W$. We have \vskip 5pt

   \begin{align*}
   p(aw)&=\frac{1}{|G|} \sum_{g \in G} \frac{1}{c(g, g^{-1})} \overline{g} p_e(\overline{g^{-1}}aw)\\
         &=\frac{1}{|G|} \sum_{g \in G} \frac{1}{c(g, g^{-1})} \overline{g} p_e((g^{-1}.a)\overline{g^{-1}}w)\\
         &=\frac{1}{|G|} \sum_{g \in G} \frac{1}{c(g, g^{-1})} \overline{g}(g^{-1}.a) p_e(\overline{g^{-1}}w)\\
         &=\frac{1}{|G|} \sum_{g \in G} \frac{1}{c(g, g^{-1})} a \overline{g}  p_e(\overline{g^{-1}}w)=ap(w).
   \end{align*}
    and
   \begin{align*}
   p(\overline{h}w)&=\frac{1}{|G|} \sum_{g \in G} \frac{1}{c(g, g^{-1})} \overline{g}p_e(\overline{g^{-1}} (\overline{h}w))\\
                   &=\frac{1}{|G|} \sum_{g \in G} \frac{c(g^{-1}, h)}{c(g, g^{-1})} \overline{g}p_e(\overline{g^{-1}h}w)\\
                   &=\frac{1}{|G|} \sum_{g \in G} \frac{c(g^{-1}, h)}{c(g, g^{-1}) c(h, h^{-1}g)} \overline{h}\hskip 2pt \overline{h^{-1}g}p_e(\overline{g^{-1}h}w)\\
                   &=\overline{h} p(w).
   \end{align*}
   where the last equality uses the following identity
  \begin{align}
   \frac{c(g^{-1}, h)}{c(g, g^{-1}) c(h, h^{-1}g)}=\frac{1}{c(h^{-1}g, g^{-1}h)}, \quad \forall g, h \in G,
  \end{align}
  which will be proved in the Appendix. So we have proved that $p$ is an $A$-module morphism. \par \vskip 5pt

  Clearly $p|_V=Id_V$. So as $A$-modules $W=V \oplus {\rm Ker}(p)$, i.e., every submodule $V$ of $W$ is a direct summand, equivalently, $W$ is a completely reducible $A$-module. Till now we have proved that every finite-dimensional $A$-module is completely reducible.\par \vskip 3pt

  Let $a=\sum_{g \in G} a_g * \overline{g}$ be a nonzero element in $A$, where each $a_g \in A_e$. Assume that $a_h \neq 0$ for some $h \in G$. Since $A_e$ is an FCR-algebra, in particular, for every nonzero element $b$ of $A_e$, there is a finite-dimensional $A_e$-module  on which $b$ acts nontrivially. So there is a finite-dimensional $A_e$-module, say $W$, such that there exists some $w \in W$ with $a_h w \neq 0$. Consider $W':=A \otimes_{A_e}W$, which is a finite-dimensioanl left $A$-module (see the notation in 4.1). Note that
  \begin{align*}
  a. (1 \otimes_{A_e}w)&=  \sum_{g \in G} a_g * \overline{g}. \hskip 2pt  1 \otimes_{A_e} w\\
                        &=\sum_{g \in G} c(g, g^{-1}) \overline{g} \otimes_{A_e} (a_gw),
  \end{align*}
  which is nonzero, since $a_h w\neq 0$ and $\{\overline{g}| g \in G\}$ is a basis of $A$ as an  $A_e$-module. Consequently, the intersection of the  annihilators of all finite-dimensional $A$-modules is zero. This completes the proof that $A$ is an FCR-algebra. \par \vskip 5pt

  For the ``only if " part, assume that $A$ is an FCR-algebra. First note that it suffices to show the case where $K$ is a splitting field of $G$. \par

  Assume that  $F$ is a finite-dimensional splitting field (see [DK]) of $G$ which contains $K$ and assume that the `` only if '' part holds when $K=F$. Write  $A_F=A \otimes_K F$ and $A_{e,F}=A_e \otimes_K F$, then $A_F$ and $A_{e, F}$ are  $F$-algebras. Moreover, we have
  \begin{align*}
  A_F=A_{e, F} *_c G.
  \end{align*}
  By [Proposition 3(2), KS] $A_F$ is an FCR-algebra over $F$. By the above assumption, $A_{e, F}$ is an FCR-algebra, and so is $A_e$, again by [Proposition 3(2), KS]. \par
  Now we can assume that $K$ is a splitting field of $G$, so $KG \simeq K \times \cdots \times K$ ($|G|$-copies). Denote  $G^*=\{\chi | \quad \chi: G \longrightarrow K^\times \mbox{ is a group morphism}\}$ to be the dual group of $G$, the order of which is $|G|$,  and $\chi(g)=1$ for all $\chi \in G^*$ if and only if $g=e$ in $G$. \par
  There is a $G^*$-action on $A$ given by
  \begin{align*}
  \chi. a * \overline{g}:= \chi(g) a * \overline{g}, \quad \forall a\in A_e, g \in G, \chi \in G^*.
  \end{align*}
  Obviously the invariant subalgebra $A^{G^*}$ is exactly $A_e$. Note that $A$ is noetherian, since it is a free $A_e$-module of finite rank and $A_e$ is noetherian. Now applying a result of Kraft and Small ([Proposition 1, KS]), we obtain that $A_e$ is an FCR-algebra. This completes the proof. \hfill $\blacksquare$
   \vskip 10pt

\sub{}\rm Let $G$ be an abelian group and $c$ a cocycle on $G$. Then the twisted group ring $KG^c=\oplus_{g \in G}Ku_g$ has multiplication given by
\begin{align*}
u_g \cdot u_h=c(g, h) u_{gh}.
\end{align*}
Note that $KG^c$ is a cocycle twist of the group algebra $KG$ (see 3.4). \par \vskip 10pt

As an application of Theorem 4.2, we have
\begin{prop}
Let $A=\oplus_{g\in G}A_g$ be a $G$-graded algebra where $G$ is a finite abelian group, and let $A^c$ be the cocycle twist of $A$ (see {\rm 3.4}). Assume that both $A$ and $A^c$ are noetherian.  Then $A$ is an FCR-algebra if and only if $A^c$ is an FCR-algebra.
\end{prop}

\noindent{\bf Proof}\quad For a homogeneous element $a$ of $A$ (or $A^c$), write $|a|$ for the degree of $a$. There is a $G$-action on $A^c$ as algebra automorphisms given by
\begin{align*}
g. a^c := \frac{c(g, |a|)}{c(|a|, g)} a^c,
\end{align*}
for all homogeneous elements $a^c \in A^c$ and $g \in G$. (Here $a^c$ denotes the corresponding element in $A^c$ of $a \in A$, see 3.4.) Note that $B_c$ (see 2.1) is bilinear, so the above action is well-defined. With this action and the cocycle $c$, we define the crossed product $A^c *_c G$.\par \vskip 5pt

Define a map $\Psi: A^c *_c G \longrightarrow KG^c \otimes A$ by
\begin{align*}
\Psi(a^c * \overline{g})=c(|a|, g) u_{|a|g} \otimes a
\end{align*}
for all homogeneous elements $a^c \in A^c$ and $g \in G$.  Clearly $\Psi$ is bijective with the inverse given by \begin{align*}
\Psi^{-1}(u_g \otimes a)= \frac{1}{c(|a|, |a|^{-1}g)} a ^c* \overline{|a|^{-1}g}.
\end{align*} \par \vskip 5pt

We claim that $\Psi$ is an algebra map, hence an isomorphism. Then if $A$ is an FCR-algebra, so is $KG^c \otimes A$ by [Proposition 3(1), KS]. Via the ismorphism $\Psi$,  apply Theorem 4.2, we get that  $A^c$ is an FCR-algebra. Conversely, note that $(A^c)^{c^{-1}}=A$, so if  $A^c$ is FCR  then $A$ is FCR, as required.\par \vskip 5pt

To see that $\Psi$ is an algebra map, note that in $A^c*_c G$,
\begin{align*}
a^c*\overline{g} \cdot b^c *\overline{h}&= c(g, h) a^c\circ (g.b^c) *\overline{gh}\\
                                    &=c(g, h)c(|a|, |b|) (a(g.b^c))* \overline{gh}\\
                                    &=c(g, h)c(|a|, |b|) \frac{c(g, |b|)}{c(|b|, g)} (ab)^c * \overline{gh},
\end{align*}
where ``$\circ$" denotes the multiplication in $A^c$ and $g.b^c$ denotes the $G$-action on $A^c$.\par

So
\begin{align*}
\Psi(a^c *\overline{g} \cdot b^c *\overline{h})= \frac{c(g, |b|)}{c(|b|, g)}c(g, h)c(|a|, |b|)c(|a||b|, gh) u_{|a||b|gh }\otimes ab.
\end{align*}
On the other hand,
\begin{align*}
\Psi(a^c*\overline{g}) \cdot \Psi(b^c* \overline{h})=c(|a|, g)c(|b|, h)c(|a|g, |b|h) u_{|a||b|gh} \otimes ab.
\end{align*}
So the fact that $\Psi$ is an algebra map, i.e., $\Psi(a^c * \overline{g}) \cdot \Psi(b^c * \overline{h})=\Psi(a^c * \overline{g} \cdot b^c * \overline{h})$, follows from the following identity
\begin{align}
{c(g, |b|)}c(g, h)c(|a|, |b|)c(|a||b|, gh)=c(|b|, g)c(|a|, g)c(|b|, h)c(|a|g, |b|h),
\end{align}
which will be proved in the Appendix. So this completes the proof. \hfill $\blacksquare$
\par \vskip 10pt

\begin{rem}
If $G$ is a finite abelian group and $c$ a cocycle of $G$,  then we have that $A$ is noetherian if and only if $A^c$ is  noetherian. \par
 Recall a fact: if $R=\oplus_{g \in G}R_g$ is a noetherian $G$-graded algebra, then $R_e$ is also noetherian. \par

In fact, if $A$ is noetherian, so is $KG^c \otimes A$. Therefore $A^c *_c G$ is
noetherian via the isomorphism $\Psi$. Note that $A^c * G=\oplus_{g \in G} A^c * \overline{g}$ is a  $G$-graded algebra with the $e$-th component $A^c$. Applying  the above fact, we get that $A^c$ is noetherian. Conversely, just note that $(A^c)^{c^{-1}}=A$, i.e., $A$ is a cocycle twist of $A^c$ by $c^{-1}$.
\end{rem}

\vskip 10pt

 \sub{} \rm Let $A=\oplus_{g \in G}A_g$ and $B=\oplus_{g \in G}B_g$ be $G$-graded algebras. Let $\varepsilon$ be a bilinear form on $G$. Then the $\varepsilon$-tensor product $A\overline{\otimes}^\varepsilon B$ of $A$ and $B$ is known as
 \begin{align*}
 A\overline{\otimes}^\varepsilon B= \oplus_{g, h \in G} A_g \otimes B_h
 \end{align*}
 with multiplication
 \begin{align*}
 a \otimes b \cdot a' \otimes b':= \varepsilon(|b|, |a'|)aa' \otimes bb'
 \end{align*}
 for  all homogeneous elements $a, a' \in A$ and $b, b' \in B$. So we have the following result, which can be viewed as an analogon of [Proposition 3(1), KS].\par \vskip 10pt

 \begin{prop}
 Assume further that $G$ is a finite abelian group and $A \otimes B$ is a noetherian algebra. If $A$ and $B$ are FCR-algebras, then $A\overline{\otimes}^\varepsilon B$ is also an FCR-algebra.
 \end{prop}

\noindent{\bf Proof} \quad Note that the usual tensor product $A \otimes B$ is a $G \times G$-graded algebra with the decomposition of vector spaces as
\begin{align*}
A \otimes B =\oplus_{(g, h) \in G \times G} A_g \otimes B_h.
\end{align*}
 Define a cocycle c on $G \times G$  as
\begin{align*}
c((g, h), (g', h')):=\varepsilon(h, g')
\end{align*}
for all $g, g', h, h' \in G$. So we have the cocycle twist $(A\otimes B)^c$ of $A\otimes B$. \par

An observation is that
\begin{align*}
A\overline{\otimes}^\varepsilon B \simeq (A\otimes B)^c,
\end{align*}
 sending $a\otimes b$ to $(a\otimes b)^c$ (where $(a \otimes b)^c$ denotes the elements of $(A \otimes B)^c$, see 3.4). Now the results follows from Proposition 4.4 and Remark 4.5.  \hfill $\blacksquare$

\vskip 10pt

\sub{Proof of Theorem 1.1}\rm \vskip 5pt
\noindent (1). Note that by Corollary  2.6 and Remark 4.5, $\mathcal{U}(L)$, $\mathcal{U}(L^c)$ and $\mathcal{U}(L)^c$ are noetherian. Now (1) follows directly  from Theorem 3.5 and Proposition 4.4.\par \vskip 5pt
\noindent (2).  There is a $G$-action on the augmented enveloping algebra $\mathcal{U}(L)$  given by
   \begin{center}
   $g. x:=\varepsilon(g, h)x, \quad \forall x \in \mathcal{U}(L)_h, g, h \in G$.
   \end{center}
   With this action and the trivial cocycle of $G$, which will be denoted by $1$, we can define the crossed product $\mathcal{U}(L)*_1 G$. \par
    Note that the following map is an algebra isomorphism:
   \begin{align*}
   \Phi: \widetilde{\mathcal{U}}(L) \longrightarrow \mathcal{U}(L)*_1 G
   \end{align*}
 where $\Phi(g \otimes x)= \varepsilon(g, |x|)x * \overline{g}$, for all homogeneous elements $x \in \mathcal{U}(L)$ and $g \in G$. (Recall that $\mathcal{U}(L)$ is a $G$-graded algebra, see 2.4, and $|x|$ denotes the degree of $x$.)  Clearly $\Phi$ is bijective, and so  it suffices to show that $\Phi$ is an algebra map. We have
 \begin{align*}
 \Phi((g \otimes x) \cdot (h \otimes y))&= \Phi (\varepsilon(|x|, h) gh \otimes xy)\\
 &= \varepsilon(|x|, h) \varepsilon(gh, |xy|) xy * \overline{gh}\\
 &=\varepsilon(g, |xy|)\varepsilon(h, |y|)xy * \overline{gh}.
 \end{align*}
And
 \begin{align*}
\Phi(g\otimes x) \cdot \Phi(h \otimes y)&= \varepsilon(g, |x|) \varepsilon(h, |y|)\hskip 2pt (x * \overline{g}) \cdot (y* \overline{h})\\
&= \varepsilon(g, |x|) \varepsilon(h, |y|)\hskip 2pt x (g.y)* \overline{gh}\\
& = \varepsilon(g, |x|) \varepsilon(h, |y|) \varepsilon(g, |y|) \hskip 2pt xy*\overline{gh},
 \end{align*}
 where $x, y \in \mathcal{U}(L)$ are homogeneous elements and $g, h \in G$. So we obtain that $\Phi((g \otimes x) \cdot (h \otimes y))=\Phi(g\otimes x) \cdot \Phi(h \otimes y)$, i.e., $\Phi$ is an algebra map.

  Now the result follows directly from Theorem 4.2. This completes the proof of Theorem 1.1.   \hfill $\blacksquare$
\vskip 20pt

\section{{\bf The color Lie algebra $sl_2^c$}}
\vskip 10pt
 In this section, we will study the finite-dimensional representations of the color Lie algebra $sl_2^c$ (see Example 3.3), which is a cocycle twist of the Lie algebra $sl_2=sl(2, K)$. The field $K$ is assumed to be algebraically closed of characteristic zero. \par \vskip 5pt

 \sub{} \rm Recall that $sl_2^c$ is a $G$-graded $\varepsilon$-Lie algebra, where $G=\mathbb{Z}_2 \times \mathbb{Z}_2$ and
 \begin{align*}
 \varepsilon((\alpha_1, \alpha_2),(\beta_1, \beta_2))=(-1)^{\alpha_1\beta_2-\alpha_2\beta_1}
 \end{align*}
 for all $\alpha_1, \alpha_2, \beta_1, \beta_2 \in \mathbb{Z}_2$. A homogeneous basis of $sl_2^c$ is given as $\{a_1, a_2, a_3\}$ such that $|a_1|=(1, 0)$, $|a_2|=(0,1)$ and $|a_3|=(1,1)$. \par \vskip 5pt

 The bracket of $sl_2^c$ is given by
  \begin{align*}
\langle a_1, a_2 \rangle =a_3, \quad \langle a_2, a_3 \rangle =a_1, \quad \langle a_3, a_1 \rangle =a_2.
\end{align*}

\vskip 10pt

As we see in Example 3.3, $sl_2^c$ is a cocycle twist of $sl_2$ with respect to the triple $(c, Id, \varepsilon)$, with the cocycle $c$ given as
 \begin{align*}
c((\alpha_1, \alpha_2), (\beta_1, \beta_2))=(-1)^{\alpha_1 \beta_2}
\end{align*}
for all $\alpha_1, \alpha_2, \beta_1, \beta_2 \in \mathbb{Z}_2$. By Theorem 1.1 (1), the algebra $\mathcal{U}(sl_2^c)$ is an FCR-algebra, in particular, every finite-dimensional $\mathcal{U}(sl_2^c)$-module is completely reducible. Now we  aim to compute all the finite-dimensional simple $\mathcal{U}(sl_2^c)$-modules (up to isomorphism). \par \vskip 10pt

First list some finite-dimensional simple $\mathcal{U}(sl_2^c)$-modules, equivalently, simple representations of $sl_2^c$ as follows: \vskip 5pt

 (1)\quad $V_\alpha^n$: \quad is an $n$-dimensional simple  $\mathcal{U}(sl_2^c)$-module  with basis $\{e_1, \cdots, e_n\}$, where $\alpha=(\alpha_1, \alpha_2)$ and $\alpha_1, \alpha_2 \in \{1, -1\}$. \par
 The $\mathcal{U}(sl_2^c)$-action is given by

 \begin{align*}
 a_1.e_j&=\alpha_1\frac{(-1)^{j-1}}{2} ((2n-j)e_{j+1}-(j-1)e_{j-1}), \quad 1 \leq j \leq n-1; \\
 a_1. e_n&=\alpha_1 \frac{(-1)^{n-1}}{2}(\alpha_2 ne_n -(n-1)e_{n-1});\\
 a_2.e_j&=-\frac{1}{2} ((2n-j)e_{j+1}+ (j-1)e_{j-1}), \quad 1 \leq j \leq n-1;\\
 a_2.e_n&=-\frac{1}{2} (\alpha_2 ne_n + (n-1)e_{n-1});\\
 a_3.e_j &=\alpha_1 \frac{(-1)^j}{2} (2n-2j+1) e_j, \quad 1 \leq j \leq n,
 \end{align*}
 where we understand $e_0=0$.
 \par \vskip 5pt

 \noindent (2)\quad $W^n$:  \quad is an $n$-dimensional simple  $\mathcal{U}(sl_2^c)$-module  with basis $\{e_1, \cdots, e_n\}$, where $n$ is odd. \par
 The $\mathcal{U}(sl_2^c)$-action is given by

 \begin{align*}
 a_1.e_j &= \frac{(-1)^{j-1}}{2} ((n-j)e_{j+1}-(j-1)e_{j-1});\\
 a_2.e_j&=-\frac{1}{2} ((n-j)e_{j+1}+ (j-1) e_{j-1});\\
 a_3. e_j&=\frac{(-1)^j}{2} (n-2j+1)e_j,
 \end{align*}
 for all $1 \leq j \leq n$ , where we understand  $e_0=e_{n+1}=0$.

 \par \vskip 10pt

 Our main result of this section is
 \begin{thm}
 The color Lie algebra $sl_2^c$ has exactly four non-isomorphic simple representations $V_{\alpha }^n$ of even dimension $n$;  five non-isomorphic simple representations $V_\alpha ^n$ and $W^n$ of odd dimension $n$.
 \end{thm}
\vskip 10pt

\sub{}\rm To show Thoerem 5.2, we need the following result.\par \vskip 5pt

\begin{prop}
  There exists an algebra embedding $\Gamma: \mathcal{U}(sl_2^c) \longrightarrow M_2(\mathcal{U}(sl_2))$  given by
 \begin{align*}
 {\footnotesize \Gamma(a_1)=\frac{i}{2} \left( \begin{array}{cc}
  0 & e-f \\
   e-f & 0\end{array} \right),
 \Gamma(a_2)=-\frac{1}{2}\left( \begin{array}{cc}
  e+f & 0 \\
   0 & -e-f\end{array} \right),
 \Gamma(a_3)=\frac{i}{2}\left( \begin{array}{cc}
  0 & h \\
   -h & 0\end{array} \right),}
   \end{align*}
   where $\{e,f, h\}$ is the standard basis of $sl_2$ as in Example 3.3, and $M_2(\mathcal{U}(sl_2))$ is the algebra of full $2\times 2$ matrices with entries from $\mathcal{U}(sl_2)$.\par
   Moreover, via the embedding $\Gamma$,   $M_2(\mathcal{U}(sl_2))$ is a free $\mathcal{U}(sl_2^c)$-module of rank four.
\end{prop}

\noindent{\bf Proof}\quad Keep notation as in 5.1. Write $KG^c=\oplus_{g \in G} Ku_g$ to be the twisted group ring. First note that  there exists an isomorphism $F: KG^c \longrightarrow M_2(K)$ defined by
\begin{align*}
 {\footnotesize u_{(0,0)}\mapsto \left( \begin{array}{cc}
  1 & 0 \\
   0 & 1\end{array} \right),
 u_{(1,0)}\mapsto \left( \begin{array}{cc}
  0 & 1 \\
   1 & 0\end{array} \right),
 u_{(0,1)}\mapsto \left( \begin{array}{cc}
  1 & 0 \\
   0 & -1\end{array} \right),
 u_{(1,1)}\mapsto  \left( \begin{array}{cc}
  0 & 1 \\
   -1 & 0\end{array} \right). }
 \end{align*}

 Consider $\Gamma$  to be the composition
 \begin{align*}
 \mathcal{U}(sl_2^c) \stackrel{\Theta}{\simeq} \mathcal{U}&(sl_2)^c \stackrel{j}{\hookrightarrow} \mathcal{U}(sl_2)^c *_c G \stackrel{\Psi}{\simeq} KG^c \otimes \mathcal{U}(sl_2) \\
 &\stackrel{F\otimes Id}{\simeq} M_2(K)\otimes \mathcal{U}(sl_2) \simeq M_2(\mathcal{U}(sl_2)),
 \end{align*}
 where $\Theta$ is the isomorphism in Theorem 3.5; $j(x)=x * \overline{e}$ for all $x\in \mathcal{U}(sl_2^c)$ ($e$ is the unit of $G$); the isomorphism $\Psi$ is defined in the proof of Proposition 4.4; and $F$ is just given as above. \par \vskip 5pt

 Note that in $\mathcal{U}(sl_2)$ (and $\mathcal{U}(sl_2)^c$), we have $a_1=\frac{i}{2} (e-f)$, $a_2=-\frac{1}{2}(e+f)$ and $a_3=\frac{i}{2} h$ ($i^2=-1$). A direct computation gives the embedding $\Gamma$. (we will omit the details.) And clearly, $M_2(\mathcal{U}(sl_2))$ is a free $\mathcal{U}(sl_2^c)$-module of rank four via the embedding $\Gamma$, since $\mathcal{U}(sl_2)^c *_c G$  is a free $\mathcal{U}(sl_2)^c$-module of rank four. This completes the proof. \hfill $\blacksquare$

\par \vskip 10pt

\sub{}\rm Recall that the unique  $n+1$-dimensional simple  representation $V(n)$ of the Lie algebra $sl_2$ (e.g., see  [Hum]), is given by
\begin{align*}
e.v_j=(n-j)v_{j+1}, \quad h.v_j=(n-2j)v_{j}, \quad f.v_j=jv_{j-1},
\end{align*}
where $\{v_0, v_1, \cdots, v_n\}$ is a basis of $V(n)$ and $0 \leq j \leq n$ (Assuming  $v_{-1}=v_{n+1}$=0). \par \vskip 5pt

Consider $K^2 \otimes V(n)$ to be a left module of $M_2(K)\otimes \mathcal{U}(sl_2) \simeq M_2(\mathcal{U}(sl_2))$, where $K^2$ is a left $M_2(K)$-module in the obvious sense.  Let $\epsilon_0, \epsilon_1$ be a standard basis of  the space $K^2$ and put
\begin{align*}
v_{k, j}=\epsilon_k \otimes v_j,
\end{align*}
where $k \in \mathbb{Z}_2$ and $0 \leq j \leq n$. So  $\{v_{k, j}| k \in \mathbb{Z}_2, 0 \leq j \leq n \}$ is a basis  of $K^2 \otimes V(n)$. \par \vskip 5pt

Using  the embedding $\Gamma$ in Proposition 5.4, we obtain that there is an induced $\mathcal{U}(sl_2^c)$-action on $K^2 \otimes V(n)$ defined as  follows:
\begin{align*}
a_1.v_{k, j}&=\frac{i}{2} ((n-j)v_{k+1, j+1}-j v_{k+1, j-1}),\\
a_2.v_{k, j}&=-\frac{1}{2} (-1)^k ((n-j)v_{k, j+1}+j v_{k, j-1}),\\
a_3. v_{k, j}&=\frac{i}{2}(-1)^k (n-2j) v_{k+1, j},
\end{align*}
for all $k \in \mathbb{Z}_2$ and $ 0 \leq j \leq n$. (Note that $k$ is in $\mathbb{Z}_2$, i.e., $1+1=0$.)
\par \vskip 10pt

We have the following  observation. \par

\begin{lem}
View $K^2 \otimes V(n)$ as an $\mathcal{U}(sl_2^c)$-module as above. Then   \vskip 5pt
\noindent {\rm (1).} \quad $K^2 \otimes V(n-1) \simeq W^n \oplus W^n$ if $n$ is odd; \\
{\rm (2).} \quad $K^2\otimes V(2n-1) \simeq V_{(1, 1)}^n \oplus V_{(1, -1)}^n \oplus V_{(-1, 1)}^n \oplus  V_{(-1, -1)}^n$ for all integer $n \geq 1$.
\end{lem}

\noindent{\bf Proof}\quad (1) Assume that $n$ is odd.  Take $\{v_{k, j}\}$ to be a basis of  $K^2 \otimes V(n-1)$ as above, and  define
\begin{align*}
 u_j=v_{0, j}+ i(-1)^{j+1} v_{1, j} \quad \mbox{and}  \quad u'_j=v_{0, j}+ i(-1)^j v_{1, j}
\end{align*}
 for all $0\leq j \leq n-1$, and let $U$ and $U'$ be the subspaces spanned by $\{u_j\}$ and $\{u'_j\}$ respectively. So we have
    \begin{align*}
     K^2 \otimes V(n-1)= U \oplus U'
     \end{align*}
  as $\mathcal{U}(sl_2^c)$-modules.\par
   Moreover, $U \simeq W^n$ by sending $u_j$ to  $e_{j+1}$, $0 \leq j \leq n-1$; $U' \simeq W^n$  by sending $u'_j$ to $e_{n-j}$, $0 \leq j \leq n-1$ (here we use that fact that $n$ is odd). So the decomposition follows. \par \vskip 5pt

(2) As in (1), we define $u_j$ and $u'_j$ for $0 \leq j \leq 2n-1$  in $K^2 \otimes V(2n-1)$, and the $\mathcal{U}(sl_2^c)$-submodules $U$ and $U'$.  Still we have  $K^2 \otimes V(2n-1)=U \oplus U'$.  \par

Set
\begin{align*}
w_{j, +}:= u_j + u_{2n-j-1} \quad \mbox{and} \quad w_{j, -}:=u_j -u_{2n-j-1}
\end{align*}
for all $0 \leq j \leq n-1$. Let $U_{+}$ and $U_{-}$ be the subspaces of $U$ spanned by $\{w_{j, +}\}$ and $\{w_{j, -}\}$ respectively. Then we  get
\begin{align*}
U=U_{+} \oplus U_{-}
\end{align*}
 as $\mathcal{U}(sl_2^c)$-modules. \par

 Moreover, $U_{+}\simeq V_{(1, 1)}^n$ by sending $w_{j, +}$ to $e_{j+1}$, $0 \leq j \leq n-1$; $U_{-} \simeq V_{(1, -1)}^n$ by sending $w_{j, -}$ to $e_{j+1}$, $0 \leq j \leq n-1$. \par
  Similarly, we can prove that $U' \simeq V_{(-1, 1)}^n \oplus V_{(-1, -1)}^n$. This completes the proof of Lemma 5.6. \hfill $\blacksquare$

\vskip 10pt

\sub{Proof of Theorem 5.2}\rm \par \vskip 5pt
  First note that all the listed $\mathcal{U}(sl_2^c)$-modules  $V_\alpha^n$ and $W^n$ in 5.1   are simple and pairwise non-isomorphic.\par
   On the other hand, let $W$ be a finite-dimensional simple  $\mathcal{U}(sl_2^c)$-module. Via the embedding $\Gamma$ in Proposition 5.4, $\mathcal{U}(sl_2^c)$ is a subalgebra of $M_2(\mathcal{U}(sl_2))$ and  $M_2(\mathcal{U}(sl_2))$ is a  right free $\mathcal{U}(sl_2^c)$-module of finite rank. Write
   \begin{align*}
   W':=M_2(\mathcal{U}(sl_2))\otimes_{\mathcal{U}(sl_2^c)}W,
   \end{align*}
   which is a finite-dimensional $M_2(\mathcal{U}(sl_2))$-module (see the notation in 4.1). \par
   Consider the map
    \begin{align*}
   t: W \longrightarrow W'
   \end{align*} such that $t(w)=1 \otimes_{\mathcal{U}(sl_2^c)} w$ for all $w \in W$. Clearly $t$ is injetive and  it is an $\mathcal{U}(sl_2^c)$-module morphism. Since $\mathcal{U}(sl_2^c)$ is an FCR-algbera (see 5.1), $W$ is a direct summand of $W'$ as $\mathcal{U}(sl_2^c)$-modules. \par \vskip 5pt

   Recall a well-known fact: let $A$ be an associative $K$-algebra and $M_2(A)$ denotes the algebra of $2 \times 2$ mtrices with entries from $A$. Then for any left $M_2(A)$-module $M$, there is an $A$-module $N$ such that $M \simeq K^2 \otimes N$. \par

  Apply this fact to the $M_2(\mathcal{U}(sl_2))$-module $W'$. So there exists a finite-dimensional $\mathcal{U}(sl_2)$-module $V$ such that $K^2 \otimes V \simeq W'$. Write $V=\oplus V(n_i)$ as the decompostion of $\mathcal{U}(sl_2)$-modules for some non-negative numbers $n_i$. So we get
  \begin{align*}
  W' \simeq \oplus K^2 \otimes V(n_i)
  \end{align*}
  as $M_2(\mathcal{U}(sl_2))$-modules. Consequently, $W$ is a direct summand of $\oplus K^2 \otimes V(n_i)$ as $\mathcal{U}(sl_2^c)$-modules. Now applying  the decompositions in Lemma 5.6, we see that $W$ is among the list of $V_\alpha^n$ and $W^n$. This completes the proof. \hfill $\blacksquare$

\par \vskip 20pt

\section{{\bf Appendix}}
\vskip 10pt

In this section we will give the detailed proofs of the identity  (4.1) and (4.2), which are direct consequences of the cocycle condition (see the equation (3.1)). However it seems that there is no exact reference. \par

Throughout $G$ will be a finite abelian group with identity elementt $e$ and $c: G \times G \longrightarrow K^\times$ a cocycle, i.e.,
\begin{align*}
c(g, h) c(gh, k)= c(h, k) c(g, hk), \quad \forall g, h, k \in G.
\end{align*}
We denote the above equation by $(g, h, k)$. \par \vskip 5pt

\sub{}\rm By $(g,e,e)$ and $(e, e, h)$, we get
\begin{align*}
c(g, e)=c(e, h)=c(e, e), \quad \forall g, h \in G.
\end{align*}

By $(h, h^{-1}, h)$, i.e.,
\begin{align*}
c(h, h^{-1})c(e, h)=c(h^{-1}, h)c(h, e),
\end{align*}
 we get $c(h, h^{-1})=c(h^{-1}, h)$ for all $h \in G$. \par \vskip 10pt

\sub{}\rm The identity (4.1) is equivalent to
\begin{align*}
c(g^{-1}, h)c(h^{-1}g,, g^{-1}h)=c(h, h^{-1}g) c(g, g^{-1}), \quad \forall g, h \in G.
\end{align*}

In fact
\begin{align*}
c(g^{-1}, h)c(h^{-1}g, g^{-1}h)&= c(h^{-1}g, g^{-1}) c(h^{-1}, h) \quad (\mbox{use } (h^{-1}g, g^{-1}, h) )\\
                                &=c(h^{-1}g, g^{-1})c(h, h^{-1}) \quad (\mbox{use } c(h, h^{-1})=c(h^{-1}, h) )\\
                                &=c(h, h^{-1}g)c(g, g^{-1}) \quad (\mbox{use } (h, h^{-1}g, g^{-1}) ).
\end{align*}
This proves the identity  (4.1). \par \vskip 10pt

\sub{} \rm For the identity (4.2),
\begin{align*}
&c(g, |b|) c(g, h) c(|a|, |b|) c(|a||b|, gh)\\
&=c(g, |b|)c(g, h) c(|b|, gh)c(|a|, |b|gh) \quad (\mbox{use } (|a|, |b|, gh))\\
&=c(g, |b|)c(|b|, g)c(|b|g, h) c(|a|, |b|gh)\quad (\mbox{use } (|b|, g, h))\\
&=c(|b|, g) c(|b|, h)c(g, |b|h) c(|a|, |b|gh)\quad (\mbox{use } |b|g=g|b| \mbox{ and } (g, |b|, h))\\
&=c(|b|, g)c(|b|, h)c(|a|, g)c(|a|g, |b|h) \quad (\mbox{use }|b|gh=g |b|h \mbox{ and }(|a|, g, |b|h)  ).
\end{align*}
This completes the proof.

\vskip 10pt

\bibliography{}

\vskip 10pt

\end{document}